
\magnification1200
\input amstex.tex
\documentstyle{amsppt}

\hsize=12.5cm
\vsize=18cm
\hoffset=1cm
\voffset=2cm

\footline={\hss{\vbox to 2cm{\vfil\hbox{\rm\folio}}}\hss}
\nopagenumbers
\def\DJ{\leavevmode\setbox0=\hbox{D}\kern0pt\rlap
{\kern.04em\raise.188\ht0\hbox{-}}D}

\def\txt#1{{\textstyle{#1}}}
\baselineskip=13pt
\def\hf{{\textstyle{1\over2}}}
\def\a{\alpha}\def\b{\beta}
\def\d{{\,\roman d}}
\def\e{\varepsilon}
\def\f{\varphi}
\def\G{\Gamma}

\def\t{\theta}
\def\={\;=\;}

\def\zt{\zeta(\hf+it)}

\def\D{\Delta}

\def\z{\zeta}

 \def\t{\theta}
\def\hf{{\textstyle{1\over2}}}
\def\txt#1{{\textstyle{#1}}}
\def\f{\varphi}

\font\tenmsb=msbm10
\font\sevenmsb=msbm7
\font\fivemsb=msbm5
\newfam\msbfam
\textfont\msbfam=\tenmsb
\scriptfont\msbfam=\sevenmsb
\scriptscriptfont\msbfam=\fivemsb
\def\Bbb#1{{\fam\msbfam #1}}

\def \NN {\Bbb N}

\def \RR {\Bbb R}

\font\ff=cmr8
\def\txt#1{{\textstyle{#1}}}
\baselineskip=13pt

\font\teneufm=eufm10
\font\seveneufm=eufm7
\font\fiveeufm=eufm5
\newfam\eufmfam
\textfont\eufmfam=\teneufm
\scriptfont\eufmfam=\seveneufm
\scriptscriptfont\eufmfam=\fiveeufm
\def\mathfrak#1{{\fam\eufmfam\relax#1}}

\font\tenmsb=msbm10
\font\sevenmsb=msbm7
\font\fivemsb=msbm5
\newfam\msbfam
     \textfont\msbfam=\tenmsb
      \scriptfont\msbfam=\sevenmsb
      \scriptscriptfont\msbfam=\fivemsb
\def\Bbb#1{{\fam\msbfam #1}}

\def \NN {\Bbb N}

\def \RR {\Bbb R}

  \def\rightheadline{{\hfil{\ff
  On the Riemann zeta-function and the divisor problem}\hfil\tenrm\folio}}

  \def\leftheadline{{\tenrm\folio\hfil{\ff
   A. Ivi\'c }\hfil}}
  \def\emptyheadline{\hfil}
  \headline{\ifnum\pageno=1 \emptyheadline\else
  \ifodd\pageno \rightheadline \else \leftheadline\fi\fi}

\topmatter
\title
ON THE RIEMANN ZETA-FUNCTION AND THE DIVISOR PROBLEM
\endtitle
\author   Aleksandar Ivi\'c  \endauthor
\address
Aleksandar Ivi\'c, Katedra Matematike RGF-a
Universiteta u Beogradu, \DJ u\v sina 7, 11000 Beograd,
Serbia (Yugoslavia).
\endaddress
\email {\tt
aivic\@rgf.bg.ac.yu,  aivic\@matf.bg.ac.yu} \endemail
\dedicatory
\enddedicatory
\abstract
{Let $\D(x)$ denote the error term in the Dirichlet
divisor problem, and $E(T)$ the error term in the asymptotic
formula for the mean square of $|\zt|$. If
$E^*(t) = E(t) - 2\pi\D^*(t/2\pi)$ with $\D^*(x) =
 -\D(x)  + 2\D(2x) - \hf\D(4x)$, then we obtain
$$
\int_0^T (E^*(t))^4\d t \;\ll_\e\; T^{16/9+\e}.
$$
We also show how our method of proof yields the bound
$$
\sum_{r=1}^R\left(\int_{t_r-G}^{t_r+G}|\zt|^{2}\d t\right)^4
\ll_\e T^{2+\e}G^{-2} + RG^4T^\e,
$$
where $T^{1/5+\e} \le G \ll T, \,
T < t_1 <\cdots< t_R\le 2T,\,t_{r+1}-t_r \ge 5G \,(r = 1,\ldots,R-1). $
}
\endabstract
\endtopmatter

\head
1. Introduction and statement of results
\endhead

Let, as usual,
$$
\D(x) \;=\; \sum_{n\le x}d(n) - x(\log x + 2\gamma - 1) -
{\txt {1\over4}},
\leqno(1.1)
$$
and
$$
E(T) \;=\;\int_0^T|\zt|^2\d t - T\left(\log\bigl({T\over2\pi}\bigr) + 2\gamma - 1
\right),\leqno(1.2)
$$
where $d(n)$ is the number of divisors of
$n, \gamma = -\G'(1) = 0.577215\ldots\,$
is Euler's constant. Thus $\D(x)$ denotes the error term in the
classical Dirichlet divisor problem, and $E(T)$ is the error term
in the mean square formula for $|\zt|$. An interesting analogy between $d(n)$ and
$|\zt|^2$ was pointed out by F.V. Atkinson [1] more than sixty
years ago. In his famous paper [2], Atkinson continued his research
and established an explicit
formula for $E(T)$ (see also the author's monographs
[7, Chapter 15] and [8, Chapter 2]).
The most significant terms in this formula, up to the factor $(-1)^n$, are
similar to those in Voronoi's formula (see [7, Chapter 3]) for $\D(x)$.
More precisely, in [13] M. Jutila showed that $E(T)$
should be actually compared to $2\pi\D^*(T/(2\pi))$, where
$$
\D^*(x) \;:=\; -\D(x)  + 2\D(2x) - \hf\D(4x).\leqno(1.3)
$$
Then  the arithmetic interpretation of $\D^*(x)$ (see T. Meurman [16]) is
$$
\hf\sum_{n\le4x}(-1)^nd(n) \;=\; x(\log x + 2\gamma - 1) + \D^*(x).
\leqno(1.4)
$$
We have the explicit, truncated formula (see e.g., [7] or [18])
$$
\D(x) = {1\over\pi\sqrt{2}}
x^{1\over4}\sum_{n\le N}d(n)n^{-{3\over4}}\cos(4\pi\sqrt{nx}
- {\txt{1\over4}}\pi) +
O_\e(x^{{1\over2}+\e}N^{-{1\over2}})\quad(2 \le N \ll x).
\leqno(1.5)
$$
One also has (see  [7, eq. (15.68)]), for $2 \le N \ll x$,
$$
\D^*(x) = {1\over\pi\sqrt{2}}x^{1\over4}
\sum_{n\le N}(-1)^nd(n)n^{-{3\over4}}
\cos(4\pi\sqrt{nx} - {\txt{1\over4}}\pi) +
O_\e(x^{{1\over2}+\e}N^{-{1\over2}}),
\leqno(1.6)
$$
which is completely analogous to (1.5).

\medskip
M. Jutila, in his works [13] and [14], investigated both the
local and global behaviour of
$$
E^*(t) \;:=\; E(t) - 2\pi\D^*\bigl({t\over2\pi}\bigr).
$$
He proved the mean square bound
$$
\int_{T-H}^{T+H}(E^*(t))^2\d t \ll_\e  HT^{1/3}\log^3T + T^{1+\e}
\quad(1 \ll H \le T), \leqno(1.7)
$$
which in particular yields
$$
\int_0^T(E^*(t))^2\d t \ll T^{4/3}\log^3T.\leqno(1.8)
$$
Here and later $\e$ denotes positive constants which are arbitrarily
small, but are not necessarily the same at each occurrence.
The bound (1.8) shows that, on the average, $E^*(t)$ is
of the order $\ll t^{1/6}\log^{3/2}t$,
while both $E(x)$ and $\D(x)$ are of the order $\asymp x^{1/4}$. This
follows from the mean square formulas (see e.g., [8])
$$
\int_0^T\D^2(x)\d x = (6\pi^2)^{-1}\sum_{n=1}^\infty d^2(n)n^{-3/2}T^{3/2}
+ O(T\log^4T),\leqno(1.9)
$$
and
$$
\int_0^T E^2(x)\d x = {\txt{2\over3}}(2\pi)^{-1/2}\sum_{n=1}^\infty d^2(n)n^{-3/2}T^{3/2}
+ O(T\log^4T).\leqno(1.10)
$$
The mean square formulas (1.9) and (1.10) also imply that the
inequalities $\a < 1/4$ and $\b < 1/4$  cannot hold, where $\a$ and
$\b$ are, respectively, the infima of the numbers $a$ and $b$ for which
the bounds
$$
\D(x) \ll x^{a},\qquad E(x) \ll x^{b}\leqno(1.11)
$$
hold. For upper bounds on $\a,\b$ see e.g., M.N. Huxley [5].
Classical conjectures are that $\a = \b = 1/4$ holds, although this
is notoriously difficult to prove. M. Jutila [13] succeeded in showing
the conditional estimates: if the conjectural $\a = 1/4$ holds, then this
implies that $\b \le 3/10$. Conversely, $\b = 1/4$ implies that
$\D^*(x) \ll_\e x^{\t+\e}$ holds with $\t \le 3/10$. Although one expects
the maximal orders of $\D(x)$ and $\D^*(x)$ to be approximately of the
same order of magnitude, this does seems difficult to establish.

\smallskip
In what concerns the formulas involving higher moments of $\D(x)$ and $E(t)$,
we refer the reader to the author's works  [6], [7] and [10] and
D.R. Heath-Brown [4]. In particular,  note that [10] contains proofs of
$$\eqalign{\cr
\int_0^TE^3(t)\d t &\;=\; 16\pi^4\int_0^{T\over2\pi}(\D^*(t))^3\d t
+ O(T^{5/3}\log^{3/2}T),\cr
\int_0^TE^4(t)\d t &\;=\; 32\pi^5\int_0^{T\over2\pi}(\D^*(t))^4\d t
+ O(T^{23/12}\log^{3/2}T).\cr}\leqno(1.12)
$$
In a recent work by P. Sargos and the author [12],
the asymptotic formulas of
K.-M. Tsang [19] for the cube and the fourth moment of $\D(x)$ were
sharpened to
$$
\int_1^X\D^3(x)\d x = BX^{7/4} + O_\e(X^{\b+\e})  \qquad(B > 0)\leqno(1.13)
$$
and
$$
\int_1^X\D^4(x)\d x = CX^2 + O_\e(X^{\gamma+\e})  \qquad(C > 0)\leqno(1.14)
$$
with $\b = {7\over5}, \gamma = {23\over12}$. This improves on the values
$\b = {47\over28}, \gamma  = {45\over23}$, obtained in [19].
Moreover, (1.13) and (1.14) remain valid if $\D(x)$ is replaced by $\D^*(x)$,
since their proofs used nothing more besides (1.5) and the bound $d(n)
\ll_\e n^\e$. Hence from (1.12) and the analogues of  (1.13)--(1.14)
for $\D^*(x)$, we infer then that
$$\eqalign{\cr
\int_0^TE^3(t)\d t &\;=\; B_1T^{7/4}
+ O(T^{5/3}\log^{3/2}T)\quad(B_1>0),\cr
\int_0^TE^4(t)\d t &\;=\; C_1T^2
+ O_\e(T^{23/12+\e})\quad(C_1>0).\cr}\leqno(1.15)
$$
The main aim of this paper is to provide an estimate for the upper
bound of the fourth moment of $E^*(t)$, which is the first non-trivial
upper bound for a higher moment of $E^*(t)$.
The result is the following

\medskip
THEOREM 1. {\it We have}
$$
\int_0^T (E^*(t))^4\d t \;\ll_\e\; T^{16/9+\e}.\leqno(1.16)
$$

\bigskip
Note that the bounds (1.8) and (1.16) do not seem to imply each other.
For the proof of (1.16) we shall need several lemmas, which will
be given in Section 2. The proof of  Theorem 1 will be given in
Section 3. Finally, in Section 4, it will be shown how the method of
proof of Theorem 1 can give a  proof of

\medskip
THEOREM 2. {\it Let} $T^{1/5+\e} \le G \ll T,\, T < t_1 <\cdots< t_R\le 2T,\,
t_{r+1}-t_r \ge 5G\, (r = 1, \cdots,\, R-1).
$ {\it Then}
$$
\sum_{r=1}^R\left(\int_{t_r-G}^{t_r+G}|\zt|^{2}\d t\right)^4
\ll_\e T^{2+\e}G^{-2} + RG^4T^\e.\leqno(1.17)
$$

\bigskip
The bound in (1.17) easily gives the well-known bound (see Section 4)
$$
\int_0^T|\zt|^{12}\d t \ll_\e T^{2+\e},\leqno(1.18)
$$
due to D.R. Heath-Brown [3] (who had $\log^{17}T$ instead
of the $T^\e$-factor). It is still essentially the sharpest result
concerning high moments of $|\zt|$. General sums of zeta-integrals over short
intervals, analogous to the one appearing in (1.17), were treated by
the author in [9].


\head
2. The necessary lemmas
\endhead

LEMMA 1 (O. Robert--P. Sargos [17]). {\it Let $k\ge 2$ be a fixed
integer and $\delta > 0$ be given.
Then the number of integers $n_1,n_2,n_3,n_4$ such that
$N < n_1,n_2,n_3,n_4 \le 2N$ and}
$$
|n_1^{1/k} + n_2^{1/k} - n_3^{1/k} - n_4^{1/k}| < \delta N^{1/k}
$$
{\it is, for any given $\e>0$,}
$$
\ll_\e N^\e(N^4\delta + N^2).\leqno(2.1)
$$

\medskip
LEMMA 2. {\it Let $1 \ll G \ll T/\log T$. Then we have}
$$
E(T) \le {2\over\sqrt{\pi}G}\int_0^{\infty} E(T+u)\,{\roman e}^{-u^2/G^2}\d u
+ O(G\log T),\leqno(2.2)
$$
{\it and}
$$
E(T) \ge {2\over\sqrt{\pi}G}\int_0^\infty E(T-u)\,{\roman e}^{-u^2/G^2}\d u
+ O(G\log T).\leqno(2.3)
$$

\medskip
{\bf Proof of Lemma 2}. The proofs of (2.2) and (2.3) are analogous, so only
the former will be treated in detail. From (1.2) we have, for $0 \le u\ll T$,
$$\eqalign{\cr
0 \le &\int_T^{T+u}|\zt|^2\d t = (T+u)\Bigl(\log\bigl({T+u\over2\pi}\bigr)
+ 2\gamma-1\Bigr) \cr&-  T\Bigl(\log\bigl({T\over2\pi}\bigr)
+ 2\gamma-1\Bigr) + E(T+u) - E(T).\cr}
$$
This gives
$$
E(T) \le E(T+u) + O(u\log T),
$$
hence
$$
\int_0^{G\log T} E(T)\,{\roman e}^{-u^2/G^2}\d u
\le \int_0^{G\log T} (E(T+u ) + O(u\log T))\,{\roman e}^{-u^2/G^2}\d u.
$$
The proof of (2.2) is completed when we extend the
integration to $[0,\infty)$ making a  small error, and recall that
$
\int_0^\infty {\roman e}^{-u^2/G^2}\d u = \hf\sqrt{\pi}G,\;
\int_0^\infty u{\roman e}^{-u^2/G^2}\d u = \hf G.
$

\medskip
LEMMA 3. {\it Let $1 \ll G \ll T$. Then we have}
$$
\D^*\bigl({T\over2\pi}\bigr) = {2\over \sqrt{\pi}G}\int_{0}^\infty
\D^*\bigl({T\over2\pi}\pm {u\over2\pi}\bigr)\,{\roman e}^{-u^2/G^2}\d u
+ O_\e(GT^\e).\leqno(2.4)
$$

\medskip
{\bf Proof of Lemma 3}. Both the cases of the $+$ and $-$ sign in (2.4) are
treated analogously. For example, we have
$$\eqalign{\cr&
\hf\sqrt{\pi}G\D^*\bigl({T\over2\pi}\bigr) - \int_{0}^\infty
\D^*\bigl({T\over2\pi}+{u\over2\pi}\bigr)\,{\roman e}^{-u^2/G^2}\d u  \cr&
= \int_{0}^\infty\Bigl(\D^*(T) -
\D^*\bigl({T\over2\pi} + {u\over2\pi}\bigr)\Bigr)\,
{\roman e}^{-u^2/G^2}\d u \cr&
= \int_{0}^{G\log T}\left(\D^*\bigl({T\over2\pi} ) -
\D^*({T\over2\pi} + {u\over2\pi}\bigr)\right)\,
{\roman e}^{-u^2/G^2}\d u + O(1) \cr&
\ll \int_{0}^{G\log T}\Bigl\{\Bigl|
\sum_{{2\over\pi}T\le n\le {2\over\pi}(T+u)}(-1)^nd(n)\Bigr|
+ O((1+|u|)\log T)\Bigr\}\d u
\ll_\e G^2T^\e,\cr}
$$
where we used (1.4) and $d(n) \ll_\e n^\e$. This establishes (2.4).

\medskip
The next lemma is F.V. Atkinson's classical explicit formula for $E(T)$ (see
[2], [7] or [8]).

\medskip
LEMMA 4. {\it Let $0 < A < A'$ be any two fixed constants
such that $AT < N < A'T$, and let $N' = N'(T) =
T/(2\pi) + N/2 - (N^2/4+ NT/(2\pi))^{1/2}$. Then }
$$
E(T) = \Sigma_1(T) + \Sigma_2(T) + O(\log^2T),\leqno(2.5)
$$
{\it where}
$$
\Sigma_1(T) = 2^{1/2}(T/(2\pi))^{1/4}\sum_{n\le N}(-1)^nd(n)n^{-3/4}
e(T,n)\cos(f(T,n)),\leqno(2.6)
$$
$$
\Sigma_2(T) = -2\sum_{n\le N'}d(n)n^{-1/2}(\log T/(2\pi n))^{-1}
\cos(T\log T/(2\pi n) - T + \pi /4),\leqno(2.7)
$$
{\it with}
$$
\eqalign{\cr&
f(T,n) = 2T{\roman {arsinh}}\,\bigl(\sqrt{\pi n/(2T})\bigr) + \sqrt{2\pi nT
+ \pi^2n^2} - \pi/4\cr&
=  -\txt{1\over4}\pi + 2\sqrt{2\pi nT} +
\txt{1\over6}\sqrt{2\pi^3}n^{3/2}T^{-1/2} + a_5n^{5/2}T^{-3/2} +
a_7n^{7/2}T^{-5/2} + \ldots\,,\cr}\leqno(2.8)
$$
$$\eqalign{\cr
e(T,n) &= (1+\pi n/(2T))^{-1/4}{\Bigl\{(2T/\pi n)^{1/2}
{\roman {arsinh}}\,(\sqrt{\pi n/(2T})\Bigr\}}^{-1}\cr&
= 1 + O(n/T)\qquad(1 \le n < T),
\cr}\leqno(2.9)
$$
{\it and $\,{\roman{arsinh}}\,x = \log(x + \sqrt{1+x^2}\,).$}

\medskip
LEMMA 5 (M. Jutila [13]). {\it For $A\in\RR$ we have}
$$
\cos\left(\sqrt{8\pi nT} +
\txt{1\over6}\sqrt{2\pi^3}n^{3/2}T^{-1/2} + A\right)
= \int_{-\infty}^\infty \a(u)\cos(\sqrt{8\pi n}(\sqrt{T} + u)
+ A)\d u,\leqno(2.10)
$$
{\it where $\a(u) \ll T^{1/6}$ for $u\not=0$,}
$$
\a(u) \ll T^{1/6}\exp(-bT^{1/4}|u|^{3/2}) \leqno(2.11)
$$
{\it for $u<0$, and}
$$
\a(u) = T^{1/8}u^{-1/4}\left(d\exp(ibT^{1/4}u^{3/2})
+ {\bar d}\exp(-ibT^{1/4}u^{3/2})\right) + O(T^{-1/8}u^{-7/4})\leqno(2.12)
$$
{\it for $u \ge T^{-1/6}$ and some constants $b\; (>0)$ and $d$.}

\head
3. The proof of  Theorem 1
\endhead
We shall prove that
$$
\int_T^{2T} (E^*(t))^4\d t \;\ll_\e\; T^{16/9+\e},\leqno(3.1)
$$
which easily implies (1.16) on replacing $T$ by $T/2, T/2^2,\ldots$ etc. and
summing all the results. Henceforth we assume that $T\le t\le 2T$,
$T^\e \le G \ll T^{5/12}$, and we begin by evaluating the integrals
$$
\int_{0}^\infty E(t\pm u){\roman e}^{-u^2/G^2}\d u\leqno(3.2)
$$
which appear in Lemma 2 (with $t$ replacing $T$),
truncating them at $u = G\log T$ with a negligible error.
A similar procedure was effected by D.R. Heath-Brown [4] and by
the author [7, Chapter 7], where the details of analogous
estimations may be found. It transpires that the contribution of
$\Sigma_2(T)$ (see (2.7)) in Atkinson's formula, as well as the
contribution of $n$ in $\Sigma_1(T)$ which satisfy $n > TG^{-2}\log T$
will be $\ll G\log T$, if we take in Lemma 4 $N = T$ for  $E(t)$ when
$T \le t \le2T$. What remains clearly corresponds to the truncated
formula (1.6) for $\D^*(x)$ with $N = TG^{-2}\log T$, or equivalently
$$
G\; =\; \sqrt{{T\over N}\log T}. \leqno(3.3)
$$
We combine now (2.2) with (2.4) with the $+$ sign (when $E(T) \ge 0$) or
(2.3) with (2.4) with the $-$ sign
(when $E(T) \le 0$), to obtain by the Cauchy-Schwarz inequality
$$
(E^*(t))^2  \ll_\e G^{-1}\int_{-G\log T}^{G\log T}{\roman e}^{-u^2/G^2}
(E^*(t+u))^2\d u + G^2T^\e,\leqno(3.4)
$$
provided that $T \le t \le 2T, T^\e \ll G \ll T^{5/12}$.
 Keeping in mind the
preceding discussion we thus have (replacing $(t+u)^{1/4}$ with $t^{1/4}$
by Taylor's formula, with the error absorbed by the last term in (3.5))
by using  (1.6), (2.5), (3.3) and (3.4),
$$\eqalign{
(E^*(t))^2  \ll_\e  &\; G^{-1}\int_{-G\log T}^{G\log T}
{\roman e}^{-u^2/G^2}
(\Sigma_3^2(X;u) + \Sigma_4^2(X,N;u)
+ \Sigma_5^2(X,N;u))\d u \cr&+ T^{1+\e}N^{-1},\cr}\leqno(3.5)
$$
where we set
$$\eqalign{\cr&
\Sigma_3(X;u) :=  t^{1/4}\times\cr&
\sum_{n\le X}
(-1)^n d(n)n^{-3/4}\Bigl\{e(t+u,n)\cos(f(t+u,n)) -
\cos(\sqrt{8\pi n(t+u)}-\pi/4)\Bigr\},\cr}\leqno(3.6)
$$
$$\eqalign{\cr&
\Sigma_4(X,N;u) := t^{1/4}\sum_{X<n\le N}
(-1)^n d(n)n^{-3/4}e(t+u,n)\cos(f(t+u,n)),\cr&
\Sigma_5(X,N;u) := t^{1/4}
\sum_{X<n\le N}(-1)^n d(n)n^{-3/4}\cos(\sqrt{8\pi n(t+u)}-\pi/4),
\cr}\leqno(3.7)
$$
where we suppose that ($N = N(T)$ is the analogue of $N$ in (1.5)
and (1.6) (cf. (3.4)), and not of $N$ in Lemma 4)
$$
T^\e \le X < T^{1/3},\; \max(X, T^{1/6}\log T) < N \ll T^{11/17}.\leqno(3.8)
$$
Here $X = X(T)$ is a parameter which allows one (by using (2.8))
to replace, in $\Sigma_3(X;u)$, $\cos(f(t+u,n))$ by
$$
(1 + cn^{3/2}(t+u)^{-1/2})\cos(\sqrt{8\pi n(t+u)}-\pi/4)
$$
plus terms of a lower order of magnitude. Note that, for $n\le X
\;(< T^{1/3})$, we may also replace $e(t+u,n)$ in
 (2.6) by 1 with the error absorbed by the
last term in (3.5). The conditions imposed in (3.8) imply that $G$
(see (3.4)) satisfies $G \ll T^{5/12}$.  Hence instead of $\Sigma_3(X;u)$
in (3.5), we may estimate
$$
\Sigma_6(X;u) := \sum_{n\le X}t^{-1/4}
(-1)^n d(n)n^{3/4} \cos(\sqrt{8\pi n(t+u)}-\pi/4),\leqno(3.9)
$$
which has the advantage because the
cosine contains $\sqrt{8\pi n(t+u)}-\pi/4$
instead of the more complicated function $f(t+u,n)$. Thus with the aid of
(3.5)--(3.9) we see that the left-hand side of (3.1) is
majorized by the maximum taken over $|u| \le G\log T$ times
$$
\eqalign{\cr&
\int_T^{2T} (E^*(t))^2(\Sigma^2_4(X,N;u) + \Sigma^2_5(X,N;u) +
\Sigma^2_6(X;u) + T^{1+\e}N^{-1})\d t\cr& \ll_\e
{\left\{\int_T^{2T} (E^*(t))^4\d t\,\int_T^{2T}\Bigl(\Sigma^4_4(X,N;u) +
\Sigma^4_5(X,N;u)+\Sigma^4_6(X;u)\Bigr)\d t\right\}}^{1/2}\cr&
 + T^{7/3+\e}N^{-1},\cr} \leqno(3.10)
$$
where we used the Cauchy-Schwarz inequality for integrals and (1.8). Thus
from (3.10) we have the key bound
$$\eqalign{
\int_T^{2T} (E^*(t))^4 \d t &\ll_\e
\max_{|u|\le G\log T}\int_T^{2T}\Bigl(\Sigma^4_4(X,N;u) +
\Sigma^4_5(X,N;u)+\Sigma^4_6(X;u)\Bigr)\d t \cr& + T^{7/3+\e}N^{-1}.
\cr}\leqno(3.11)
$$
To evaluate the integrals on the right-hand side of (3.11)
 we note first that
$$
\int_T^{2T}\Bigl(\Sigma^4_4(X,N;u) + \ldots\Bigr)\d t
\le \int_{T/2}^{5T/2}\f(t)\Bigl(\Sigma^4_4(X,N;u) + \ldots\Bigr)\d t,
\leqno(3.12)
$$
where $\f(t)$ is a smooth, nonnegative function supported in
$\,[T/2,\,5T/2]\,$, such that $\f(t) = 1$ when $T \le t \le 2T$.
The integrals of $\sum^4_4(X,N;u), \sum^4_5(X,N;u)$ and
$\sum^4_6(X;u)$ are all estimated analogously.
The sums over $n$ are divided into
$\ll \log T$ subsums of the form $\sum_{K<n\le K'\le 2K}$, the cosines
are written as exponentials, and the fourth power is written as a quadruple
sum over the integer variables $m,n,k,l$. Then we perform
a large number of integrations by parts to deduce that the contribution
of those $m,n,k,l$ for which $|\D| \ge T^{\e-1/2}$ is negligible
(i.e., $\ll T^{-A}$ for any fixed $A>0$), where
$$
\D \;:=\; \sqrt{8\pi}(\sqrt{m} + \sqrt{n} - \sqrt{k} - \sqrt{l}\,).
\leqno(3.13)
$$
Therefore, in the case of $\Sigma_5(X,N;u)$, there  remains the estimate
$$
\eqalign{\cr&
\int_T^{2T}\Sigma^4_5(X,N;u)\d t\cr&
\ll_\e 1 + T^{1+\e}\max_{|u|\le G\log T}\sup_{X\le K \le N}
{\int_{T/2}^{5T/2}}\f(t)\times\cr&
\Bigl|{\mathop{\sum\nolimits^*}\limits_{K<m,n,k,l\le K'\le 2K}}
(-1)^{m+n+k+l}d(m)d(n)d(k)d(l)(mnkl)^{-3/4}\exp(i\D\sqrt{t+u})\Bigr|\d t,
\cr}\leqno(3.14)
$$
where $\sum^* $ means that $|\D| \le T^{\e-1/2}$ holds. Now we use
Lemma 1 (with $k=2$, $\delta \asymp K^{-1/2}|\D|$), estimating the integral
on the right-hand side of (3.14) trivially.
We obtain that the left-hand side of (3.14) is
$$
\eqalign{\cr&
\ll_\e T^{1+\e}\max_{X\le K\le N,|\D|\le T^{\e-1/2}}
K^{-3}T(K^4K^{-1/2}|\D| + K^2)\cr&
\ll_\e T^\e(T^2N^{1/2}T^{-1/2} + T^2X^{-1})\cr&
\ll_\e T^{3/2+\e}N^{1/2} + T^{2+\e}X^{-1}.\cr}\leqno(3.15)
$$
Proceeding  analogously as in (3.15), we obtain that
$$\eqalign{\cr
\int_T^{2T}\Sigma^4_6(X;u)\d t &\ll_\e
T^{1+\e}\max_{1\le K\le X,|\D|\le T^{\e-1/2}}T^{-1}K^3(K^4K^{-1/2}|\D|
+ K^2)\cr&
\ll_\e T^\e(T^{-1/2}X^{13/2} + X^5),\cr}\leqno(3.16)
$$
since instead of $(mnkl)^{-3/4}$ in (3.14) now we shall have
$(mnkl)^{3/4}t^{-1}$  (see (3.9)).

The estimation of $\Sigma_4(X,N;u)$ (see (3.7))
presents a technical problem, since the
cosines contain the function $f(t,n)$, and Lemma 1 cannot be applied
directly. First we note that, by using (2.8), we can expand the exponential
in power series to get rid of the terms $a_5n^{5/2}t^{-3/2} + \ldots\;$.
 In this
process the main term will be 1, and the error terms will make a contribution
which will be (for shortness we set
$a = \sqrt{8\pi}, b = \txt{1\over6}\sqrt{2\pi^3}$ and $\tau = t+u$)
$$\eqalign{\cr&
\ll_\e \max_{|u|\le G\log T}\sup_{X\le K\le N}T\int_{T/2}^{5T/2}\f(t)
\Bigl|\sum_{K<n\le2K}(-1)^nd(n)n^{7/4}\tau^{-3/2}\times\cr&\times
\exp\Bigl(ia(n\tau)^{1/2} +
 ib(n^{3}/\tau)^{1/2}\Bigr)\Bigr|^4\d t\cr&
\ll_\e \max_{|u|\le G\log T}\sup_{X\le K\le N}
T^{\e-5}K^{9/2}\int_{T/2}^{5T/2}\f(t)
\Bigl|\sum_{K<n\le2K}(-1)^nd(n)n^{7/4}
\times\cr&\times\exp\Bigl(ia(n\tau)^{1/2} +
 ib(n^{3}/\tau)^{1/2}\Bigr)\Bigr|^2\d t\cr&
\ll_\e \max_{|u|\le G\log T}\sup_{X\le K\le N}
T^{\e-5}K^{9/2}(T\sum_{K<n\le2K}n^{7/2} \cr&
\quad\quad+ T^{1/2}\sum_{K<m\not=n\le2K}
(mn)^{7/4}|\sqrt{m}-\sqrt{n}|^{-1})\cr&
\ll_\e \max_{|u|\le G\log T}\sup_{X\le K\le N}
T^{\e-5}K^{9/2}(TK^{9/2} +
+ T^{1/2}K^4\sum_{K<m\not=n\le2K}|m-n|^{-1})\cr&
\ll_\e \max_{X\le K\le N}T^{\e-5}K^{9/2}TK^{9/2} \ll_\e T^{\e-4}N^9 \ll_\e
T^{3/2+\e}N^{1/2}
\cr}
$$
for
$
N \;\ll\; T^{11/17},
$
which is implied by (3.8). Thus we are left with
$$
\cos\left(\sqrt{8\pi n\tau} +
\txt{1\over6}\sqrt{2\pi^3}n^{3/2}\tau^{-1/2}-\txt{1\over4}\pi\right)
$$
in  $\Sigma_4(X,N;u)$, and we can apply Lemma 5. With $\a(v)$ given by
(2.12) we have
$$
\eqalign{\cr&
\cos\left(\sqrt{8\pi n\tau} +
\txt{1\over6}\sqrt{2\pi^3}n^{3/2}\tau^{-1/2}- A\right)
= O(T^{-10})\, +
\cr&
 \int_{-u_0}^{u_1}\a(v)\cos(\sqrt{8\pi n}(\sqrt{\tau}+v) -
A)\d v  +
\int_{u_1}^\infty \a(v)\cos(\sqrt{8\pi n}(\sqrt{\tau}+v) -
A)\d v ,\cr}\leqno(3.17)
$$
where we set
$$
u_0 =  T^{-1/6}\log T, \; u_1 = CKT^{-1/2},\leqno(3.18)
$$
and $C>0$ is a large constant.

\smallskip We proceed now as in the case of $\Sigma_5(X,N;u)$. We write
the cosines as exponentials in the quadruple sum over $m,n,k,l$.
Again, after we  first perform a
large number of integrations by parts over $t$, only the portion
of the sum for which $|\D| \le T^{\e-1/2}$ will remain, where $\D$
is given by (3.13). In the
remaining sum we use (3.17) (once with $A= {1\over4}\pi$ and once with
$A = {3\over4}\pi$), noting that ${\roman e}^{iz}
= \cos z + i\cos(z - \hf\pi)$. We
remark that, for $|v| \le u_0$, we can use the crude estimate
 $\a(u) \ll T^{1/6}$, hence
for this portion the estimation will be quite analogous to
the preceding case. Next we note that
$$
\int_{u_0}^{u_1}\tau^{1/8}v^{-1/4}\exp(ib\tau^{1/4}v^{3/2}
\pm \sqrt{8\pi n}v)\d v \ll \log T\quad(\tau = t + u, \,|u| \le G\log T),
$$
writing the integral as a sum of $\ll \log T$ integrals over $[U,\,U']$
with $u_0 \le U < U' \le 2U \ll u_1$, and applying the second derivative
test to each of these integrals. We also remark
that the contribution of the $O$-term in (2.12) will be,
by trivial estimation,
$$
\int_{u_0}^{\infty}T^{-1/8}u^{-7/4}\d u  \ll T^{-1/8}u_0^{-3/4} \ll 1
$$
if we suppose that  (3.18)
is satisfied. It remains yet to
deal with the integral with $v > u_1$ in (3.17), when we note that
$$
{\partial \over \partial v}
\left(b\tau^{1/4}v^{3/2} \pm \sqrt{8\pi n}v\right) \;\gg\; T^{1/4}v^{1/2}
\quad(v > u_1),
$$
provided that $C$ in (3.18) is sufficiently large.
Hence by the first derivative test
$$
\eqalign{\cr&
\int_{u_1}^\infty \a(v)\cos(\sqrt{8\pi n}(\sqrt{\tau}+v) -
\txt{1\over4}\pi)\d v \cr&
\ll 1 + T^{1/8}u_1^{-1/4}T^{-1/4}u_1^{-1/2}
\cr& \ll 1 + T^{1/4}K^{-3/4} \ll 1 + T^{1/4}X^{-3/4} \ll 1,
\cr}
$$
since $K \gg X \gg T^{1/3}$. Thus the contribution of the integrals
on the right-hand side of (3.17) is $\ll \log T$.

Then we can proceed with
the estimation as in the case of $\Sigma_5(X,N;u)$ to obtain
$$
\int_T^{2T}\Sigma^4_4(X,N;u)\d t \ll_\e T^{3/2+\e}N^{1/2} + T^{2+\e}X^{-1}.
$$
Gathering together all the bounds, we see that the integral in (3.1)
is
$$
\ll_\e T^\e\Bigl(T^{3/2}N^{1/2} + T^2X^{-1}
+ T^{-1/2}X^{13/2} + X^5
 + T^{7/3+\e}N^{-1}\Bigr),\leqno(3.19)
$$
provided that (3.8) holds. Finally we choose
$$
X \;=\; T^{1/3-\e},
\quad N \;=\; T^{5/9},
$$
so that (3.8) is fulfilled. The above terms are then $\ll_\e T^{16/9+\e}$,
and the proof of Theorem 1 is complete. The limit of the method is the
bound $\ll T^2X^{-1} \ll T^{5/3}$, which would yield the exponent
$5/3+\e$ in (1.16). The true order of the integral in (1.16), and in general
the order of the $k$-th moment of $E^*(t)$, is elusive. This comes
from the definition $E^*(t) = E(t) - 2\pi\D^*(t/(2\pi))$, which makes
it difficult to see how much the oscillations of the functions $E$
and $\D^*$ cancel each other.

\head
4. The proof of Theorem 2
\endhead
We shall show now how the method of proof of our Theorem 1 may be used
to yield Theorem 2. Our starting point is an expression for the integral
$$
\int_{t_r-2G}^{t_r+2G}  \f_r(t)|\zt|^2\d t,\leqno(4.1)
$$
where $t_r$ is as in the formulation of Theorem 2, and
$\f_r \in C^\infty$ is a non-negative function supported in
$[t_r-2G,\,t_r+2G]$ that equals unity in $[t_r-G,\,t_r+G]$. The integral
in (4.1) majorizes the integral
$$
\int_{t_r-G}^{t_r+G} |\zt|^2\d t,\leqno(4.2)
$$
which is of great importance in zeta-function theory (see K. Matsumoto
[15] for an extensive account on mean square theory involving $\z(s)$).
One can treat the integral in (4.1) by any of the following methods.

\hskip1cm
a) Using exponential averaging (or some other smoothing  like $\f_r$ above),
namely the Gaussian weight $\exp(-\hf x^2)$, in connection with
the function $E(T)$,
in view of F.V. Atkinson's well-known  explicit formula (cf. Lemma 4).
This is the approach employed originally by D.R. Heath-Brown [3].

\hskip1cm
b) One can use the Voronoi summation formula (e.g., see [8, Chapter 3]) for
the explicit expression (approximate functional equation) for $|\zt|^2
= \chi^{-1}(\hf + it)\zeta^2(\hf + it)$, where
 $\z(s) = \chi(s)\zeta(1-s)$, namely
$$
\chi(s) \;=\;2^s\pi^{s-1}\sin(\hf\pi s)\G(1-s).
$$
Voronoi's formula is present indirectly in Atkinson's formula, so that
this approach is more direct. The effect of the smoothing function $\f_r$
in (4.2) is to shorten the sum approximating $|\z|^2$ to the range
${T\over2\pi}(1 - G^{-1}T^\e) \le n \le {T\over2\pi} \,(T = t_r)$.
After this no integration is needed, and proceeding as in [7, Chapters 7-8]
one obtains that the integral in (4.2) equals
$O_\e(GT^\e)$ plus a multiple of
$$
\int_{t_r-2G}^{t_r+2G}  \f_r(t)
\sum_{k\le T^{1+\e}G^{-2}}(-1)^kd(k)k^{-1/2}\left({1\over4}
+ {t\over2\pi k}\right)^{-1/4}\sin f(t,k)\d t,\leqno(4.3)
$$
where $f(t,k)$ is given by (2.8).

\hskip1cm
c) Instead of the Voronoi summation formula one can use the (simpler)
Poisson summation formula, namely
$$
\sum_{n=1}^\infty f(n) = \int_0^\infty f(x)\d x +
2\sum_{n=1}^\infty  \int_0^\infty f(x)\cos(2\pi n x)\d x,
$$
provided that $f(x)$ is smooth and compactly supported in $(0,\,\infty)$.
In [11] a sketch of this approach is given.

\medskip
We begin now the derivation of (1.17), simplifying first in (4.3) the factor
$(1/4 + t/(2\pi k))^{-1/4}$ by Taylor's formula, and then
raising the expression in (4.3) to the fourth power, using H\"older's
inequality for integrals. It follows that the sum in (1.17) is bounded by
$$\eqalign{\cr&
RG^4T^\e + T^{-1}G^3\sum_{r=1}^R\int\limits_{t_r-2G}^{t_r+2G}  \f_r(t)
\Bigl|\sum_{k\le T^{1+\e}G^{-2}}(-1)^kd(k)k^{-1/4}\sin f(t,k)\Bigr|^4\d t
\cr&
\ll_\e RG^4T^\e + T^{-1}G^3\int_{T/2}^{5T/2}\f(t)
\Bigl|\sum_{k\le T^{1+\e}G^{-2}}(-1)^kd(k)k^{-1/4}\sin f(t,k)\Bigr|^4\d t,
\cr}\leqno(4.4)
$$
where $\f(t)$
is a non-negative, smooth function supported in $[T/2,5T/2]$ such that
$\f(t) = 1$ for $T \le t \le 2T$, hence $\f^{(m)}(t) \ll_m T^{-m}$.
Therefore it suffices to bound the expression
$$
I_K := \int_{T/2}^{5T/2}\f(t)
\Bigl|\sum_{K<k\le K'\le2K}(-1)^kd(k)k^{-1/4}{\roman e}^{if(t,k)}
\Bigr|^4\d t,\leqno(4.5)
$$
where $T^{1/3} \le K \ll T^{1+\e}G^{-2}$, $T^{1/5+\e} \le G \le T^{1/3}$.
Namely for $K \le T^{1/3}$ the contribution is trivially
$\ll RG^4T^\e$, and the same holds (e.g., see [7, Theorem 7.3])
for the values $G \ge T^{1/3}$. Recall that
$$
f(t,k) = -\txt{1\over4}\pi + 2\sqrt{2\pi kt} +
\txt{1\over6}\sqrt{2\pi^3}k^{3/2}t^{-1/2} + a_5k^{5/2}t^{-3/2} +
a_7k^{7/2}t^{-5/2} + \ldots\,,
$$
and note that we have $k^{5/2}t^{-3/2}
\ll T^{1+\e}G^{-5} \le T^{-\e}$ for $G \ge T^{1/5+\e}$. This means that
we may replace, on the right-hand side of (4.5), $f(t,k)$ in the exponential
by
$$
-\txt{1\over4}\pi + 2\sqrt{2\pi kt} +
\txt{1\over6}\sqrt{2\pi^3}k^{3/2}t^{-1/2}
$$
times a series whose terms are of descending order of magnitude.
The main contribution will thus come from the above term.

\bigskip
After this procedure we see that
the integral in (4.5) bears close resemblance to the integral of the
fourth moment of $E^*(t)$. The term $k^{3/2}t^{-1/2}$ in the exponential
is treated by the use of Lemma 5, similarly as was done in the case of
$\Sigma_4(X,N;u)$ in Section 3. In our case, due to the fact
that $K \ge T^{1/3}$ may be assumed, there will be no sum corresponding
to $\Sigma_3(X;u)$. Now we proceed similarly as in the proof of  Theorem 1.
We shall apply  Lemma 5 as in the proof of  Theorem 1.
Developing  the fourth power in (4.5) and performing a large number
of integrations by parts, we see that only the values for which
$$
|E| \;\le\; T^{\e-1/2},\quad
E = \sqrt{8\pi}(\sqrt{m} + \sqrt{n} - \sqrt{k} - \sqrt{l}\,)
$$
will be relevant, where $m,n,k,l$ are integers from $\,[K,\,K']$.
Thus, by Lemma 1 (with $\delta = T^{-1/2+\e}K^{-1/2}$) and
trivial estimation, their contribution to $I_K$ will be
$$\eqalign{\cr&
\ll_\e T^{1+\e}K^{-1}(K^4T^{-1/2}K^{-1/2} + K^2)\cr&
\ll_\e T^{1+\e}K^{5/2}T^{-1/2} \ll_\e T^{3+\e}G^{-5}.
\cr}
$$
This yields the bound
$$
\sum_{r=1}^R\left(\int_{t_r-G}^{t_r+G}|\zt|^{2}\d t\right)^4
\ll_\e RG^4T^\e + G^3T^{-1}T^{3+\e}G^{-5},
$$
which is (1.17).

\medskip
It remains to show how (1.17) gives the twelfth moment estimate
(1.18). Write
$$
\int_T^{2T}|\zt|^{12}\d t \le \sum_{r\le T+1}|\z(\hf + i\tau_r^*)|^{12},
  \leqno(4.6)
$$
where for $r = 1,2,\ldots$ we set
$$
|\z(\hf + i\tau_r^*)| \;:=\; \max_{T+r-1\le t \le T+r} |\zt|.
$$
Let $\{t_{r,V}\}$ be the subset of $\{\tau_r^*\}$ such that
$$
V \le |\z(\hf + it_{r,V})| \le 2V\qquad(r = 1,\ldots,\,R_V),
$$
where clearly $V$ may be restricted to $O(\log T)$ values of the form
$2^m$ such that $\log T\le V \le T^{1/6}$, since $\zt = o(t^{1/6})$
(see [7, Chapter 7]). Now since we have (see e.g., [8, Theorem 1.2]),
for fixed $k\in\NN$,
$$
|\zt|^k \ll \log t\int_{t-{1\over2}}^{t+{1\over2}}|\z(\hf + it + iu)|^k\d u
+ 1,
$$
it follows that, for some points $t_r'\,(\in\,[T,2T])$
with $r = 1,\ldots\,,R',\,
R' \le R_V$, $1 \ll G \ll T, \,t'_{r+1}-t_r' \ge 5G$,
$$\eqalign{\cr&
R_VV^{2} \le \sum_{r=1}^{R_V}|\z(\hf + it_{r,V})|^{2}\cr&
\ll \sum_{r=1}^{R_V}\log T\Bigl(\int_{t_{r,V}-{1\over2}}^{t_{r,V}+{1\over2}}
|\zt|^{2}\d t + 1\Bigr)\cr&
\ll \sum_{r=1}^{R'}\log T\Bigl(\int_{t_r'-G}^{t_r'+G}|\zt|^{2}\d t\Bigr)
+ R_V\log T\cr&
\le \log T(R')^{3/4}\Bigl(
\sum_{r=1}^{R'}\int_{t_r'-G}^{t_r'+G}|\zt|^{2}\d t\Bigr)^{1/4}
+ R_V\log T\cr&
\ll_\e T^\e(R_VG + R_V^{3/4}T^{1/2}G^{-1/2}),
\cr}
$$
where the estimate of Theorem 2 was used, with $R_V$ replacing $R$.
If we take $G = V^2T^{-2\e}$, then we obtain
$$
R_V^{1/4} \ll_\e T^{1/2+\e}G^{-3/2},
$$
which gives
$$
 R_V \ll_\e T^{2+\e}G^{-6} \ll_\e T^{2+\e}V^{-12}.
$$
Then the portion of the sum in (4.6) for which $|\z(\hf+i\tau_r*)| \ge
T^{1/10+\e}$ is
$$
\ll \log T\max_{V \ge T^{1/10+\e}}R_VV^{12} \ll_\e T^{2+\e},
$$
But for values of $V$ such that $V \le T^{1/10+\e}$, the above
bound easily follows from the large values estimate (the fourth moment)
 $R \ll_\e T^{1+\e}V^{-4}$.
This shows that the integral in (4.6) is $\ll_\e T^{2+\e}$, and
proves (1.18). Note that the author [9, Corollary 1] proved the bound
$$
\sum_{r=1}^R\Bigl(\sum_{t_r-G}^{t_r+G}|\zt|^4\d t\Bigr)^2
\ll RG^2\log^8T + T^2G^{-1}\log^CT\leqno(4.7)
$$
for some $C>0$, where $T < t_1 <\ldots < t_R \le 2T\,$,
$t_{r+1} - t_r \ge 5G$ for $r = 1,\ldots , R-1$ and
$1 \ll G \ll T$. The bound (4.7), which
is independent of Theorem 2, was proved by a method different from
the one used in this work. Like (1.17), the bound (4.7) also leads
to the twelfth moment estimate (1.18).

\vfill
\eject
\topglue2cm
\Refs
\bigskip

\item{[1]} F.V. Atkinson, The mean value of the zeta-function on the
critical line, Quart. J. Math. Oxford {\bf10}(1939), 122-128.

\item{[2]} F.V. Atkinson, The mean value of the Riemann zeta-function,
Acta Math. {\bf81}(1949), 353-376.

\item{[3]} D.R. Heath-Brown, The twelfth power moment of the Riemann
zeta-function, Quart. J. Math. (Oxford) {\bf29}(1978), 443-462,

\item{[4]} D.R. Heath-Brown, The distribution of moments
in the Dirichlet divisor problems, Acta Arith. {\bf60}(1992),
389-415.

\item{[5]} M.N. Huxley, Area, Lattice Points and Exponential
Sums, Oxford Science Publications, Clarendon Press,
Oxford, 1996

\item{[6]} A. Ivi\'c, Large values of the error term in the
divisor problem, Invent. Math. {\bf71}(1983), 513-520.

\item{[7]} A. Ivi\'c, The Riemann zeta-function, John Wiley \&
Sons, New York, 1985.

\item{[8]} A. Ivi\'c, The mean values of the Riemann zeta-function,
LNs {\bf 82}, Tata Inst. of Fundamental Research, Bombay (distr. by
Springer Verlag, Berlin etc.), 1991.

\item{[9]} A. Ivi\'c, Power moments of the Riemann zeta-function
over short intervals, Arch. Mat. {\bf62}\ (1994),\ 418-424.

\item{[10]} A. Ivi\'c, On some problems involving the mean square
of $|\zt|$, Bull. CXVI Acad. Serbe 1998, Classe des Sciences
math\'ematiques {\bf23}, 71-76.

\item{[11]} A. Ivi\'c, Sums of squares of $|\zt|$ over short
intervals, Max-Planck-Institut f\"ur Mathematik, Preprint
Series 2002({\bf52}), 12 pp.

\item{[12]} A. Ivi\'c and P. Sargos, On the higher moments of the
error term in the divisor problem, to appear.

\item{[13]} M. Jutila, Riemann's zeta-function and the divisor problem,
Arkiv Mat. {\bf21}(1983), 75-96 and II, ibid. {\bf31}(1993), 61-70.

\item{[14]} M. Jutila, On a formula of Atkinson, in ``Proc. Coll. Soc.
J. Bolyai'' Vol. {\bf34} (Budapest, 1981), North-Holland,
Amsterdam, 1984, 807-823.

\item{[15]} K. Matsumoto, Recent developments in the mean square theory
of the Riemann zeta and other zeta-functions, in ``Number Theory",
Birkh\"auser, Basel, 2000, 241-286.

\item{[16]} T. Meurman, A generalization of Atkinson's formula to
$L$-functions, Acta Arith. {\bf47}(1986), 351-370.

\item{[17]} O. Robert and P. Sargos, Three-dimensional
exponential sums with monomials, J. reine angew. Math. (in print).

\item{[18]} E.C. Titchmarsh, The theory of the Riemann zeta-function
(2nd ed.),  University Press, Oxford, 1986.

\item{[19]} K.-M. Tsang, Higher power moments of $\D(x)$, $E(t)$
and $P(x)$, Proc. London Math. Soc. (3){\bf 65}(1992), 65-84.

\bigskip
\leftline{Aleksandar Ivi\'c}
\leftline{Katedra Matematike RGF-a}
\leftline{Universiteta u Beogradu, \DJ u\v sina 7,}
\leftline{11000 Beograd, Serbia (Yugoslavia)}
\leftline{e-mail: {\tt aivic\@matf.bg.ac.yu, aivic\@rgf.bg.ac.yu}}

\endRefs


\bye